\documentclass[11pt]{amsart}

\usepackage[T1]{fontenc}
\usepackage[utf8]{inputenc}
\usepackage{lmodern}
\usepackage{microtype}

\usepackage{amsmath,amssymb,amsfonts}
\usepackage{mathtools}
\usepackage{xparse}
\usepackage{mathrsfs}

\usepackage[authoryear,round]{natbib}

\theoremstyle{plain}
\newtheorem{theorem}{Theorem}[section]
\newtheorem{proposition}[theorem]{Proposition}
\newtheorem{lemma}[theorem]{Lemma}
\newtheorem{corollary}[theorem]{Corollary}

\theoremstyle{definition}
\newtheorem{definition}[theorem]{Definition}
\newtheorem{example}[theorem]{Example}

\theoremstyle{remark}
\newtheorem{remark}[theorem]{Remark}

\newcommand{\cD}{\mathcal{D}}

\newcommand{\cH}{\mathcal{H}}

\newcommand{\abs}[1]{\lvert #1\rvert}

\NewDocumentCommand{\norm}{m g}{%
	\lVert #1\rVert\IfValueT{#2}{_{#2}}%
}
\NewDocumentCommand{\Norm}{m g}{%
	\left\lVert #1\right\rVert\IfValueT{#2}{_{#2}}%
}

\renewcommand{\d}{\mathop{}\!\mathrm{d}}

\numberwithin{equation}{section}

\begin{document}
	
	\title[Spectral Criteria for Uniqueness Pairs]
	{Spectral Criteria for Uniqueness Pairs of Unitary Transforms}
	
	\author{Oleg Szehr}
	\address{Dalle Molle Institute for Artificial Intelligence (IDSIA), SUPSI/USI,
		Via la Santa 1, 6962 Lugano-Viganello, Switzerland}
	\email{oleg.szehr@idsia.ch}
	
	\begin{abstract}
		The identification of sampling sets that enable unique signal recovery is fundamental
		to many applications in signal processing and remains a central problem in mathematical
		analysis. Recent studies, particularly in the context of the Fourier transform and
		crystalline measures, have developed a theory of recovery from two-sided sampling,
		where samples are prescribed simultaneously in the physical and transformed domains.
		Kulikov, Nazarov, and Sodin introduced a method for identifying such uniqueness pairs
		based on functional inequalities of the Wirtinger-Poincar\'e type. In this work, we
		propose an alternative spectral approach motivated by quantum mechanics. The guiding
		observation is that zeros of a function and of its transform impose Dirichlet-type
		confinement in two conjugate representations, thereby converting two-sided uniqueness
		questions into lower-bound problems for confined Hamiltonians. This reveals a confinement
		form of the uncertainty principle: excessive simultaneous confinement in two conjugate
		representations forces the underlying state to vanish. For the Fourier transform, the
		relevant Hamiltonian is the harmonic oscillator, and the uniformly supercritical
		uniqueness criterion is recovered through a variational spectral argument. Our viewpoint
		extends to unitary transforms whose associated localization operators admit local
		Sturm-Liouville or Schr\"odinger-type confined realizations, a class that includes
		transforms commonly used in signal processing and mathematical physics. It abstracts
		the Wirtinger-Poincar\'e mechanism by replacing the ordinary Dirichlet-Laplacian
		constant with the local spectral floor of a Hamiltonian-type operator associated with
		the transform. We formulate this principle for Sturm--Liouville operators with weights
		or nontrivial coefficients, and illustrate it for the fractional Fourier transform and
		the Hankel transform, where phase-space rotation and singular endpoint behavior enter
		the uniqueness criteria.
	\end{abstract}
	
	
	
	\maketitle

\section{Introduction}\label{sec1}
The ability to reconstruct a function from its values on a discrete set lies at the heart of sampling theory, a cornerstone of both digital signal processing and applied mathematics. A central role in sampling theory is played by the concept of a \emph{uniqueness set}, which generalizes the classical Nyquist rate of Whittaker-Kotelnikov-Shannon sampling theory beyond Paley-Wiener spaces. The characterization of uniqueness sets is closely tied to the geometry of the sampling space and to the density of sampling points~\citep{S}.

Recently, particularly in the context of the Fourier transform and crystalline measures, the concept of a \emph{uniqueness pair} has emerged. Such a pair consists of two sets, one in the physical domain and one in the Fourier domain, with the property that prescribed data on these sets determine the function uniquely; see~\citep{RV}, and also~\citep{RS,KNS}. Beyond its mathematical significance, this concept has practical implications in signal recovery, since it addresses uniqueness from \emph{two-sided} sampling, where samples are taken simultaneously in the original and transformed domains. Achieving uniqueness typically requires that the sampling densities in both domains exceed certain critical thresholds, whose precise form depends on the geometry of the underlying spaces.

Current approaches to determining these recovery thresholds rely crucially on estimates derived from the Wirtinger-Poincar\'e inequality~\citep{KNS}. Motivated by analogies with quantum mechanics, the present work develops a complementary spectral viewpoint. The guiding interpretation is that zeros of $f$ impose Dirichlet boundary conditions in the physical variable, while zeros of $\widehat f$ impose Dirichlet boundary conditions in the frequency variable. Thus a two-sided uniqueness problem can be interpreted as a confinement problem for a quantum Hamiltonian. For the Fourier transform, the natural Hamiltonian is the harmonic oscillator
$$
H=\frac12\left(-\frac{d^2}{dx^2}+x^2\right),
$$
whose quadratic-form domain is the Fourier-symmetric Sobolev space considered below. The intervals determined by the zero sets decompose the relevant quadratic forms into confined Dirichlet problems, and uniqueness follows when the corresponding local Dirichlet spectral floors are uniformly supercritical. This spectral formulation recovers the Fourier uniqueness criterion from a variational principle rather than directly from the Wirtinger-Poincar\'e inequality. It also gives the critical density threshold a physical interpretation: simultaneous confinement above the threshold in two conjugate representations forces the underlying quantum state to vanish.

Beyond the classical Fourier transform, several unitary transforms arising in quantum mechanics, signal processing, and mathematical physics admit analogous two-sided sampling questions. We formulate a spectral approach to such questions under an explicit locality assumption: the localization operator associated with the transform must allow Dirichlet confinement on the intervals determined by the imposed zeros. This assumption is satisfied in important cases where the confined operators have Sturm-Liouville or Schr\"{o}dinger-type realizations. In these settings, lower bounds on the respective local spectral floors provide sufficient conditions for uniqueness pairs.

This shows that the spectral mechanism is not merely the Wirtinger-Poincar\'e argument in different notation. In the Fourier case, the Wirtinger-Poincar\'e constant appears because the relevant local spectral floor is that of the confined Dirichlet Laplacian. In more general settings, the local spectral floor is determined by the full confined operator associated with the transform. This distinction is important in cases where the confinement energy also reflects potential terms, endpoint behavior, or Sturm-Liouville structure.

To further illustrate the scope of the approach, we consider the fractional Fourier and Hankel transforms. The fractional Fourier transform provides a continuous interpolation between the time and frequency domains and is widely used in optical signal analysis, radar imaging, and time-frequency filtering~\citep{OZK}. The Hankel transform arises naturally in problems with radial symmetry, including wave propagation, diffraction, and medical imaging such as computed tomography~\citep{KS}.

\section{Uniqueness pairs for the Fourier transform}\label{sec:ordinaryFourier}

We begin by considering the Fourier-symmetric Sobolev space $$\cH=\left\{f\:\Big|\: \int_{\mathbb{R}}(1+x^2)\abs{f(x)}^2\d x+ \int_{\mathbb{R}}(1+k^2)\abs{\hat{f}(k)}^2\d k<\infty\right\},$$ which arises naturally in the theory of \textit{modulation spaces} - function spaces used in time-frequency analysis to quantify the joint localization of a signal in both time and frequency; see, e.g.,~\citep{EG}. More recently, $\mathcal{H}$ has also attracted attention in the study of crystalline measures and two-sided sampling theory~\citep{KNS,Z}. In this section, we examine $\mathcal{H}$ from a quantum-mechanical perspective and recover several recent results concerning uniqueness pairs associated with the Fourier transform. We adopt the unitary convention for the Fourier transform (defined initially on Schwartz test functions and then through unitary extension to $L^2(\mathbb R)$), as it facilitates the generalization of the notion of uniqueness pairs to other unitary transformations:
\begin{align*}
	\hat{f}(k) &= {F}[f](k) = \frac{1}{\sqrt{2\pi}}\int_{\mathbb{R}}f(x)e^{-ikx}\d x,\\
	f(x) &= {F}^{-1}[\hat{f}](x)=\frac{1}{\sqrt{2\pi}}\int_{\mathbb{R}}\hat{f}(k)e^{ikx}\d k.
\end{align*}
We note that the convention used in~\citep{KNS} is not unitary. It will be convenient to view $\mathcal{H}$ as a subspace of $L^2(\mathbb{R})$, by introducing the norm 
$$\Norm{f}{\mathcal{H}}^2=\int_{\mathbb{R}}x^2\abs{f(x)}^2\d x+ \int_{\mathbb{R}}k^2\abs{\hat{f}(k)}^2\d k$$ 
and interpreting $\cH=\{f\in L^2(\mathbb{R})\:|\:\norm{f}{\mathcal{H}}<\infty\}$. Defining
$$H=\frac{1}{2}\left(-\frac{d^2}{d x^2}+x^2\right)$$
one computes formally using Plancherel's identity and integration by parts
\begin{align*}
	\Norm{f}{\mathcal{H}}^2&=\int_{\mathbb{R}}x^2\abs{f(x)}^2\d x+ \int_{\mathbb{R}}k^2\abs{\hat{f}(k)}^2\d k=\int_{\mathbb{R}}x^2\abs{f(x)}^2\d x+ \int_{\mathbb{R}}\abs{f'(x)}^2\d x\\
	&=\int_{\mathbb{R}}\overline{f(x)}x^2f(x)\d x- \int_{\mathbb{R}}\overline{f(x)}\frac{\partial^2f}{\partial x^2}(x)\d x=2\int_{\mathbb{R}}\bar{f} H(f)\d x.
\end{align*}
$H$ is the Hamiltonian operator of the classical quantum harmonic oscillator~\citep{NS}. The eigenvalue-eigenvector equation $Hf=Ef$ is the Schr\"{o}dinger equation of the system.

The preceding computation is justified for compactly supported smooth functions. On that set the operator $H$ is essentially self-adjoint and positive; we denote its self-adjoint closure again by \(H\). Since \(k\widehat f(k)\in L^2\) is equivalent, by
Plancherel, to \(f'\in L^2\) in the weak sense, closure gives
\[
\mathcal H=\mathcal Q(H)=\mathcal D(H^{1/2}),
\]
where  \(\mathcal D(H)\) denotes the
operator domain and \(\mathcal Q(H)\) is the quadratic-form domain of the closed form
\[
 h[f]
=
\frac12\int_{\mathbb R}
\bigl(|f'(x)|^2+x^2|f(x)|^2\bigr)\,\d x,
\qquad
\operatorname{dom} h=\mathcal H.
\]
Equivalently, $\|f\|_{\mathcal H}^2=2 h[f]=2\|H^{1/2}f\|^2_{{L^2}(\mathbb R)}$. The described passage between a positive self-adjoint operator to its associated closed
quadratic form and form domain is standard; see, for instance,~\cite[Chapter VI]{Kato}.

\begin{theorem}\label{thm:standardHamiltonian}
	The following assertions hold for the Fourier-symmetric Sobolev space $\mathcal{H}$:
	\begin{enumerate}
		\item \(\mathcal H\subset L^2(\mathbb R)\) and \(\mathcal H\) is both the
		domain of \(H^{1/2}\) and the image of \(H^{-1/2}\), where $H=\frac12\left(-\frac{d^2}{dx^2}+x^2\right)$
		is the positive self-adjoint harmonic oscillator.
		\item $H$ has eigenfunctions that are Hermite functions of the form $$\phi_n(x)=\frac{1}{\sqrt{n!}}\left(\frac{1}{\sqrt{2}}(x-\partial_x)\right)^n\phi_0(x),\ n\geq0,$$ with $\phi_0=\pi^{-1/4}e^{-x^2/2}$. The respective eigenvalues are $E_n=(n+\frac{1}{2})$.
		\item $H$ is self-dual under the Fourier transform, in the sense that $\hat{H}:={F}H{F}^{-1}=H$ as self-adjoint operators on $L^2(\mathbb R)$.
		\item The Hermite functions form an orthonormal basis of \(L^2(\mathbb R)\)
		and an orthogonal basis of \(\mathcal H\) with respect to the form inner
		product.
		\item For all $f\in\mathcal{H}$, $f\neq0$, it holds that $\norm{f}{\mathcal{H}}\geq\norm{f}{L^2(\mathbb{R})}$ and equality is achieved if and only if $f$ is a multiple of $\phi_0$.
	\end{enumerate}
\end{theorem}
The properties \emph{(1), (2), (3)} are classical facts in the study of the quantum oscillator, see e.g.~\citep{Kato,NS} and see \citep{PS} for a nice summary. Property \emph{(4)} is known in the theory of modulation spaces, cf.~\citep{EG}; see also \citep{J1} for further elaboration, and has recently appeared in the setting of $\mathcal{H}$ in~\citep{Z}. In the present framework, \emph{(4)} follows from \emph{(1)} and
\emph{(2)}: since \(H^{-1/2}:L^2(\mathbb R)\to\mathcal H\) is an
isomorphism and \(H^{-1/2}\phi_n=E_n^{-1/2}\phi_n\), the Hermite basis of
\(L^2(\mathbb R)\) is carried to a rescaling of itself. Property~\emph{(5)} is an uncertainty principle in \(\mathcal H\): the two
integrals in the definition of \(\|f\|_{\mathcal H}\) cannot be
simultaneously small. It follows from the spectral estimate
\[
\|f\|_{\mathcal H}^2
=
2 h[f]
=
2\|H^{1/2}f\|_{L^2(\mathbb R)}^2
\geq
2E_0\|f\|_{L^2(\mathbb R)}^2
=
\|f\|_{L^2(\mathbb R)}^2,
\]
where \(E_0=1/2\). Equality holds if and only if \(f\) belongs to the
ground-state eigenspace, which is one-dimensional and spanned by
\(\phi_0\).

A central theme of our analysis is the spectral behavior of Hamiltonians
confined by infinite potential barriers. For a bounded interval \([a,b]\),
we write \(H_{[a,b]}\) for the self-adjoint and positive operator on
\(L^2([a,b])\) acting by
\begin{align*}
H_{[a,b]}f
=
\frac12\left(-\frac{d^2f}{dx^2}+x^2f\right)
\end{align*}
under Dirichlet boundary conditions $f(a)=f(b)=0$. As for the unconfined oscillator, this construction can be made rigorous
by means of the associated closed quadratic form
\[
 h_{[a,b]}[f]
=
\frac12\int_a^b\bigl(|f'(x)|^2+x^2|f(x)|^2\bigr)\,\d x,
\]
defined on the space of square-integrable functions on \((a,b)\) with
square-integrable weak derivative and vanishing boundary values
\(f(a)=f(b)=0\).

We shall impose pointwise zero conditions on functions in \(\mathcal H\) and
on their Fourier transforms. These conditions are meaningful: for every
\(f\in\mathcal H\), Plancherel gives \(f'\in L^2(\mathbb R)\) in the weak
sense, and similarly \((\widehat f)'\in L^2(\mathbb R)\). Hence \(f\) and
\(\widehat f\) have continuous representatives. Consequently, if
\(f(a)=0\) and \(f(b)=0\), then the restriction of \(f\) to
\([a,b]\) belongs to the form domain of the confined
Hamiltonian \(H_{[a,b]}\). The same observation applies to
\(\widehat f\). In this sense,
the imposed zeros act as Dirichlet walls.

Bearing this construction in mind we formally write
\begin{align}
	Hf=Ef\ \text{under}\ f(a)=f(b)=0\label{confinedSchroed}
\end{align}
for the Schr\"{o}dinger equation with Dirichlet boundary conditions.

The following definition follows~\citep{KNS}. Item~\emph{1.} is taken
from there\footnote{The notions of super- and subcriticality studied
	in~\citep{KNS} are more general in that we restrict ourselves to the
	symmetric situation \((p,q)=(2,2)\), cf.~Def.~2 in~\citep{KNS}. This is
	to achieve a lean presentation.}, whereas \emph{2.} and \emph{3.} are original to this
work.
\begin{definition}\label{criticalDef}
	Let \(\Lambda=\{\lambda_j:j\in\mathbb Z\}\) and
	\(M=\{\mu_j:j\in\mathbb Z\}\) be discrete subsets of \(\mathbb R\)
	with strictly increasing enumerations that decompose \(\mathbb R\) into
	the intervals \((\lambda_j,\lambda_{j+1})\) and
	\((\mu_j,\mu_{j+1})\), respectively.
	\begin{enumerate}
		\item The pair $(\Lambda,M)$ is uniformly supercritical if
		\begin{align*}\sup_{j\in\mathbb{Z}}\max\{\abs{\lambda_j},\abs{\lambda_{j+1}}\}(\lambda_{j+1}-\lambda_j)&<\pi\ \textnormal{and}\\ \sup_{j\in\mathbb{Z}}\max\{\abs{\mu_j},\abs{\mu_{j+1}}\}(\mu_{j+1}-\mu_j)&<\pi.\end{align*}
		The pair $(\Lambda,M)$ is supercritical if
		\begin{align*}\limsup\abs{\lambda_j}(\lambda_{j+1}-\lambda_j)&<\pi\ \textnormal{and}\\ \limsup\abs{\mu_j}(\mu_{j+1}-\mu_j)&<\pi.\end{align*}
		The pair $(\Lambda,M)$ is subcritical if
		\begin{align*}\liminf\abs{\lambda_j}(\lambda_{j+1}-\lambda_j)&>\pi\ \textnormal{and}\\
			\liminf{\abs{\mu_j}}(\mu_{j+1}-\mu_j)&>\pi.\end{align*}
		Here, $\limsup$ and $\liminf$ are meant as two-sided conditions.
		\item Let \(E_{[a,b]}^{(0)}\) denote the ground-state energy of the harmonic
		oscillator confined to a bounded interval \([a,b]\) with Dirichlet boundary
		conditions, cf.~equation~\eqref{confinedSchroed}. The pair $(\Lambda,M)$ has uniformly supercritical ground-state energy if 
		$$ \inf_{j\in\mathbb{Z}}\frac{\sqrt{E_{[\lambda_j,\lambda_{j+1}]}^{(0)}}}{\max\{\abs{\lambda_j},\abs{\lambda_{j+1}}\}}>1\quad \textnormal{and}\quad\inf_{j\in\mathbb{Z}}\frac{\sqrt{E_{[\mu_j,\mu_{j+1}]}^{(0)}}}{\max\{\abs{\mu_j},\abs{\mu_{j+1}}\}}>1.$$
		The corresponding supercritical and subcritical ground-state energy
		conditions are defined analogously, cf.~point \emph{1.}~above.
		
		\item For a bounded interval \([a,b]\subset\mathbb R\), let \(\nu_{[a,b]}\) denote the
		lowest eigenvalue, in the variational sense, of the
		Sturm-Liouville problem with Dirichlet boundary conditions
		\[
		-\frac{d^2f}{dx^2}=\nu f,\qquad f|_{\partial [a,b]}=0 .
		\]
		Equivalently,
		\[
		\nu_{[a,b]}
		=
		\inf_{f\neq0}
		\frac{\int_{[a,b]} |f'(x)|^2\,\d x}
		{\int_{[a,b]} |f(x)|^2\,\d x},
		\]
		where the infimum is over the space of square-integrable functions on $I$ with square-integrable weak derivative and vanishing boundary conditions. The pair \((\Lambda,M)\) has uniformly supercritical
		Sturm-Liouville spectrum if
		$$\inf_{j\in\mathbb{Z}}\frac{\sqrt{\nu_{[\lambda_j,\lambda_{j+1}]}}}{\max\{\abs{\lambda_j},\abs{\lambda_{j+1}}\}}>1\quad \textnormal{and}\quad\inf_{j\in\mathbb{Z}}\frac{\sqrt{\nu_{[\mu_j,\mu_{j+1}]}}}{\max\{\abs{\mu_j},\abs{\mu_{j+1}}\}}>1.$$
		The corresponding supercritical and subcritical Sturm-Liouville
		conditions are defined analogously.
	\end{enumerate}
\end{definition}

\begin{remark}\label{rem:scalingPairs}
	A pair $(\tilde{\Lambda},\tilde{M})$ is called uniformly supercritical in~\citep{KNS} if 
	\begin{align*}
\sup_{j\in\mathbb{Z}}\max\{\abs{\tilde\lambda_j},\abs{\tilde\lambda_{j+1}}\}(\tilde\lambda_{j+1}-\tilde\lambda_j)&<\frac{1}{2},\\ \sup_{j\in\mathbb{Z}}\max\{\abs{\tilde\mu_j},\abs{\tilde\mu_{j+1}}\}(\tilde\mu_{j+1}-\tilde\mu_j)&<\frac{1}{2}.
\end{align*}
	This discrepancy is due to the choice of Fourier convention. To translate between the conventions consider $M'=2\pi\tilde M$ such that $\hat{f}(\mu_i')=0$ for all $\mu_j'\in M'$ under our convention. The respective criticality condition for $M'$ reads $\sup_{j\in\mathbb{Z}}\max\{\abs{\mu_j'},\abs{\mu_{j+1}'}\}(\mu_{j+1}'-\tilde\mu_j')<\frac{1}{2}\cdot(2\pi)^2$. $(\tilde{\Lambda},M')$ is a uniqueness pair if and only if $(\alpha\tilde\Lambda,\frac{1}{\alpha} M')$, $\alpha\neq0$ is a uniqueness pair. We use this property setting $\Lambda=\sqrt{2\pi}\tilde{\Lambda}$ and $M=\frac{1}{\sqrt{2\pi}}M'$.
\end{remark}

The three notions of Definition~\ref{criticalDef} are closely related. The
Sturm-Liouville formulation (item \emph{3.}) is simply the spectral version of the
criticality condition of~\cite{KNS} (item \emph{1.}), while the ground-state condition is
a closely related condition for the confined harmonic oscillator (item \emph{2.}). The precise
comparison is as follows.

\begin{lemma}\label{lem:criticalComparison}
	Let \((\Lambda,M)\) be as in Definition~\ref{criticalDef}. Then the following
	statements hold.
	\begin{enumerate}
		\item Uniform supercriticality is equivalent to uniformly supercritical
		Sturm-Liouville spectrum.
		
		\item If \((\Lambda,M)\) has uniformly supercritical ground-state energy, then
		\((\Lambda,M)\) has uniformly supercritical Sturm-Liouville spectrum, and
		therefore is uniformly supercritical.
		
		\item If \((\Lambda,M)\) is uniformly supercritical, then it has uniformly
		supercritical ground-state energy outside at most finitely many intervals in
		each of the two sequences.
		
		\item The analogous implications hold for the nonuniform supercritical and
		subcritical conditions, with the corresponding limit superior and limit
		inferior formulations.
	\end{enumerate}
\end{lemma}
Before proving the lemma we notice that either condition implies that $(\Lambda,M)$ is a uniqueness pair.
\begin{theorem}\label{thm:mainTheorem}
	Let \((\Lambda,M)\) be as in Definition~\ref{criticalDef}. Suppose that at
	least one of the following conditions holds:
	\begin{enumerate}
		\item \((\Lambda,M)\) is uniformly supercritical;
		\item \((\Lambda,M)\) has uniformly supercritical ground-state energy;
		\item \((\Lambda,M)\) has uniformly supercritical Sturm-Liouville
		spectrum.
	\end{enumerate}
	Then the only function \(f\in\mathcal H\) satisfying $f(\lambda_i)=0$,
	$\widehat f(\mu_i)=0$, $i\in\mathbb Z$ is \(f=0\).
\end{theorem}
The first conclusion of the theorem, item~\emph{1.}, has been proved in~\citep{KNS}. Their elegant proof relies crucially on the Wirtinger-Poincar\'{e} inequality. A closely related argument, involving the same inequality, has appeared in the context of irregular sampling in~\citep{G}. We will give an alternative spectral proof and generalization of this result. But rather than making use of the Wirtinger-Poincar\'{e} inequality, we will employ spectral methods to arrive at our estimates. Our main line of reasoning is as follows: Let \([a,b]\subset\mathbb R\) and consider the operator
\[
\left(-\frac{d^2}{dx^2}\right)_{[a,b]}
\]
with Dirichlet boundary conditions on \(L^2([a,b])\). This operator is positive and self-adjoint with compact resolvent. Hence
its spectrum is discrete and consists of eigenvalues of finite multiplicity.
The corresponding Sturm-Liouville problem
\[
-f''=\nu f,\qquad f|_{\partial[a,b]}=0,
\]
has normalized eigenfunctions
\[
\psi_{[a,b]}^{(n)}(x)
=
\sqrt{\frac{2}{b-a}}
\sin\left(\frac{\pi(n+1)(x-a)}{b-a}\right),
\qquad n\geq0,
\]
with eigenvalues
\[
\nu_{[a,b]}^{(n)}
=
\frac{(n+1)^2\pi^2}{(b-a)^2}.
\]
By the variational characterization of the lowest eigenvalue,
\begin{align}
\inf_{f\neq0}
\frac{\int_{[a,b]} \bar{f}\left(-\frac{d^2}{dx^2}\right)_{[a,b]}f\,\d x}
{\int_{[a,b]} |f|^2\,\d x}=\nu_{[a,b]}^{(0)}.\label{eq:varCharLowEig}
\end{align}
Similarly, the Dirichlet realization of the confined harmonic oscillator
$H_{[a,b]}$ on \(L^2([a,b])\) is positive and self-adjoint with compact resolvent. We
denote its lowest eigenvalue by \(E_{[a,b]}^{(0)}\). The variational
characterization is
$$\inf_{f\neq0}\frac{\int_{a}^{b}\bar{f}H_{[a,b]}f\d x}{\int_{a}^{b}\abs{f}^2\d x}= E_{[a,b]}^{(0)}.$$
We now present proofs of Lemma~\ref{lem:criticalComparison} and Theorem~\ref{thm:mainTheorem}.

\begin{proof}[Proof of Lemma~\ref{lem:criticalComparison}]
	It suffices to prove the assertions for \(\Lambda\); the same argument applies
	to \(M\). By equation~\eqref{eq:varCharLowEig}
	$
	\nu_{[a,b]}=\nu_{[a,b]}^{(0)}
	=
	\frac{\pi^2}{(b-a)^2}$.
	Therefore, setting \(a=\lambda_j\) and \(b=\lambda_{j+1}\),
	\[
	\frac{\sqrt{\nu_{[\lambda_j,\lambda_{j+1}]}}}
	{\max\{|\lambda_j|,|\lambda_{j+1}|\}}>1
	\quad\Longleftrightarrow\quad
	\max\{|\lambda_j|,|\lambda_{j+1}|\}
	(\lambda_{j+1}-\lambda_j)<\pi .
	\]
	This proves the equivalence between items \emph{1.} and \emph{3.} The corresponding nonuniform
	statements follow by taking the appropriate limit superior or limit inferior.
	
	We next compare with the confined harmonic oscillator. Let
	\(E_{[a,b]}^{(0)}\) denote the Dirichlet ground-state energy of equation~\eqref{confinedSchroed}. Since \(x^2\leq \max\{|a|^2,|b|^2\}\) on \([a,b]\), the
	variational principle gives
	\[
	E_{[a,b]}^{(0)}
	\leq
	\frac12\nu_{[a,b]}+\frac12\max\{|a|^2,|b|^2\}.
	\]
	Hence, if
	\[
	\frac{E_{[a,b]}^{(0)}}{\max\{|a|^2,|b|^2\}}>1,
	\]
	then
	\[
	1<
	\frac{E_{[a,b]}^{(0)}}{\max\{|a|^2,|b|^2\}}
	\leq
	\frac12
	\frac{\nu_{[a,b]}}{\max\{|a|^2,|b|^2\}}
	+\frac12,
	\]
	and therefore
	\[
	\frac{\sqrt{\nu_{[a,b]}}}{\max\{|a|,|b|\}}>1.
	\]
	Applied to \([a,b]=[\lambda_j,\lambda_{j+1}]\), this proves that uniformly
	supercritical ground-state energy implies uniformly supercritical
	Sturm-Liouville spectrum.
	
	Conversely, suppose that
	\[
	\max\{|\lambda_j|,|\lambda_{j+1}|\}
	(\lambda_{j+1}-\lambda_j)\leq \pi-\eta
	\]
	for some \(\eta>0\) and all \(j\). For all sufficiently large \(|j|\), the
	interval \([\lambda_j,\lambda_{j+1}]\) does not contain the origin. Hence
	\(x^2\geq \min\{|\lambda_j|^2,|\lambda_{j+1}|^2\}\) on
	\([\lambda_j,\lambda_{j+1}]\), and the variational principle gives
	\[
	E_{[\lambda_j,\lambda_{j+1}]}^{(0)}
	\geq
	\frac12\nu_{[\lambda_j,\lambda_{j+1}]}
	+
	\frac12\min\{|\lambda_j|^2,|\lambda_{j+1}|^2\}.
	\]
	Hence
	\begin{align*}
	&\frac{E_{[\lambda_j,\lambda_{j+1}]}^{(0)}}
	{\max\{|\lambda_j|^2,|\lambda_{j+1}|^2\}}\\
	&\hspace{1.5cm}\geq
	\frac12
	\frac{\pi^2}
	{(\lambda_{j+1}-\lambda_j)^2
		\max\{|\lambda_j|^2,|\lambda_{j+1}|^2\}}
	+
	\frac12
	\frac{\min\{|\lambda_j|^2,|\lambda_{j+1}|^2\}}
	{\max\{|\lambda_j|^2,|\lambda_{j+1}|^2\}} .
	\end{align*}
	The first term is bounded below by $\frac12\frac{\pi^2}{(\pi-\eta)^2}$.
	Moreover, since
	\[
	\max\{|\lambda_j|,|\lambda_{j+1}|\}
	(\lambda_{j+1}-\lambda_j)\leq \pi-\eta,
	\]
	we have
	\[
	\frac{\lambda_{j+1}-\lambda_j}
	{\max\{|\lambda_j|,|\lambda_{j+1}|\}}\to0
	\qquad \text{as } |j|\to\infty,
	\]
	and therefore
	\[
	\frac{\min\{|\lambda_j|,|\lambda_{j+1}|\}}
	{\max\{|\lambda_j|,|\lambda_{j+1}|\}}\to1.
	\]
	It follows that
	\[
	\frac{E_{[\lambda_j,\lambda_{j+1}]}^{(0)}}
	{\max\{|\lambda_j|^2,|\lambda_{j+1}|^2\}}>1
	\]
	for all sufficiently large \(|j|\). Thus uniform supercriticality implies
	uniformly supercritical ground-state energy, with at most finitely many
	exceptional intervals. This proves the asserted finite-exception statement.
\end{proof}

\begin{proof}[Proof of Theorem~\ref{thm:mainTheorem}]
	Let \(f\in\mathcal H\) satisfy \(f|_\Lambda=0\) and \(\widehat f|_M=0\).
	By Lemma~\ref{lem:criticalComparison}, it suffices to prove the theorem under
	the assumption that \((\Lambda,M)\) has uniformly supercritical
	Sturm-Liouville spectrum. Hence there exists \(l>1\) such that
	\[
	\frac{\sqrt{\nu_{[\lambda_j,\lambda_{j+1}]}}}
	{\max\{|\lambda_j|,|\lambda_{j+1}|\}}\geq l
	\quad\textnormal{and}\quad
	\frac{\sqrt{\nu_{[\mu_j,\mu_{j+1}]}}}
	{\max\{|\mu_j|,|\mu_{j+1}|\}}\geq l
	\]
	for all \(j\in\mathbb Z\). By Theorem~\ref{thm:standardHamiltonian} \(\mathcal H=\mathcal D(H^{1/2})\) and
	$\|f\|_{\mathcal H}^2
	=
	2\|H^{1/2}f\|_{L^2(\mathbb R)}^2$
	and, since the Fourier transform is unitary and \(FHF^{-1}=H\), the same identity
	holds in Fourier variables $\|f\|_{\mathcal H}^2=2\|H^{1/2}\widehat f\|_{L^2(\mathbb R)}^2$. 	Averaging these two identities gives
	\[
	\|f\|_{\mathcal H}^2
	=
	\|H^{1/2}f\|_{L^2(\mathbb R)}^2
	+
	\|H^{1/2}\widehat f\|_{L^2(\mathbb R)}^2 .
	\]

We decompose \(\mathbb R\) into individual intervals that are determined by
the zeros of \(f\) and \(\widehat f\). Using the variational characterization
of the lowest Dirichlet Sturm-Liouville eigenvalue, equation~\eqref{eq:varCharLowEig},
together with the estimates
\(x^2\leq \max\{|\lambda_j|^2,|\lambda_{j+1}|^2\}\) on
\([\lambda_j,\lambda_{j+1}]\), respectively
\(k^2\leq \max\{|\mu_j|^2,|\mu_{j+1}|^2\}\) on
\([\mu_j,\mu_{j+1}]\), we obtain
\begin{align*}
	\norm{f}{\mathcal{H}}^2
	&=
	\|H^{1/2}f\|_{L^2(\mathbb R)}^2
	+
	\|H^{1/2}\widehat{f}\|_{L^2(\mathbb R)}^2\\
	&=
	\sum_{j\in\mathbb{Z}}\left(
	\|H_{[\lambda_j,\lambda_{j+1}]}^{1/2}f\|_{L^2([\lambda_j,\lambda_{j+1}])}^2
	+
	\|H_{[\mu_j,\mu_{j+1}]}^{1/2}\widehat{f}\|_{L^2([\mu_j,\mu_{j+1}])}^2
	\right)\\
	&\geq \frac{1}{2}\sum_{j\in\mathbb{Z}}\left(
	\nu_{[\lambda_j,\lambda_{j+1}]}
	\int_{\lambda_j}^{\lambda_{j+1}}\abs{f}^2\,\d x
	+
	\int_{\lambda_j}^{\lambda_{j+1}}x^2\abs{f}^2\,\d x\right.\\
	&\hspace{3cm}\left.
	+
	\nu_{[\mu_j,\mu_{j+1}]}
	\int_{\mu_j}^{\mu_{j+1}}\abs{\widehat{f}}^2\,\d k
	+
	\int_{\mu_j}^{\mu_{j+1}}k^2\abs{\widehat{f}}^2\,\d k
	\right)\\
	&\geq \frac{1}{2}\sum_{j\in\mathbb{Z}}\left(
	\left(
	\frac{\nu_{[\lambda_j,\lambda_{j+1}]}}
	{\max\{\abs{\lambda_{j}}^2,\abs{\lambda_{j+1}}^2\}}
	+1
	\right)
	\int_{\lambda_j}^{\lambda_{j+1}}x^2\abs{f}^2\,\d x \right.\\
	&\hspace{3cm}\left.
	+
	\left(
	\frac{\nu_{[\mu_j,\mu_{j+1}]}}
	{\max\{\abs{\mu_{j}}^2,\abs{\mu_{j+1}}^2\}}
	+1
	\right)
	\int_{\mu_{j}}^{\mu_{j+1}}k^2\abs{\widehat{f}}^2\,\d k
	\right).
\end{align*}
	Therefore
	\[
	\begin{aligned}
		\|f\|_{\mathcal H}^2
		\geq
		\frac{1}{2}(l^2+1)
		\left(
		\int_{\mathbb R}x^2|f(x)|^2\d x
		+
		\int_{\mathbb R}k^2|\widehat f(k)|^2\d k
		\right)
		=
		\frac{1}{2}(l^2+1)\|f\|_{\mathcal H}^2 .
	\end{aligned}
	\]
	Since \(l>1\), this is a contradiction unless \(f=0\).
\end{proof}
We emphasize some features of Theorem~\ref{thm:mainTheorem} which are not
apparent from the formulation of~\cite{KNS}. First, the criterion admits a direct quantum-mechanical interpretation. The imposed zeros of $f$ and $\widehat f$ act as infinite potential walls in position and momentum space, respectively. They therefore decompose the harmonic oscillator into confined oscillators on the intervals determined by $\Lambda$ and $M$. The harmonic-oscillator formulation in terms of the ground-state energy $E_{[a,b]}^{(0)}$ makes explicit the physical mechanism behind the estimate: sufficiently narrow confinement, relative to the distance from the origin, forces too much energy. In other words, uniqueness follows when the corresponding confinement energies force the oscillator energy to exceed what is compatible with the global $\mathcal H$-norm. 

This gives a form of uncertainty principle that is close in spirit to the classical Paley-Wiener principle. One usual formulation says that a nonzero quantum state cannot be simultaneously confined to a bounded region in position space and to a bounded region in momentum space. In other words, exact localization with respect to two incompatible spectral projections is impossible. In the setting at hand, the projections are not merely the two characteristic projections associated with two boxes. Instead, the zeros of $f$ and $\widehat f$ impose families of Dirichlet walls, and hence decompose position and momentum space into many confined regions. This leads to a finer spectral obstruction: if these local confinements are too strong, then the combined Dirichlet ground-state energies force the oscillator energy above the value allowed by the global $\mathcal H$-energy. Thus Theorem~\ref{thm:mainTheorem} may be viewed as a projection-type uncertainty principle for the Fourier transform, with the ordinary two-box localization principle replaced by simultaneous confinement relative to two families of interval projections.

Second, the argument identifies the local mechanism behind the critical threshold. On each interval, the relevant lower bound is the lowest Dirichlet Rayleigh quotient, whose extremizer is the ground state of the corresponding confined problem. Thus the spectral proof shows that the threshold is governed by the first confined eigenmode. This does not by itself prove global sharpness, but it indicates where sharpness phenomena should originate, namely from configurations in which the intervalwise ground-state bounds are approached asymptotically. 

Finally, the operator-theoretic approach suggests a generalized principle behind uniqueness pairs beyond the Wirtinger-Poincar\'{e} inequality. Its basic ingredients are the following: one needs a compatible positive self-adjoint operator, a transformation that acts on its quadratic-form domain, and a class of confinement projections for which the corresponding Dirichlet spectra can be controlled. We discuss this in the next section.
\section{Uniqueness Pairs for Unitary Transformations beyond the Fourier Transform}
\label{sec:blueprint}

We now formulate a conditional generalization of Theorem~\ref{thm:mainTheorem}. We isolate the structural assumptions behind the preceding argument and obtain a spectral criterion for uniqueness pairs associated with other unitary transforms. The two assumptions below describe the global operator-theoretic structure: a transform-invariant form domain and a localization Hamiltonian built from multiplication in the original and transformed variables. Our spectral argument also requires a separate locality condition, stated in Definition~\ref{def:localConfinement}: the relevant kinetic form must be compatible with Dirichlet confinement on the intervals determined by the imposed zeros.

Let $S\subset\mathbb R$ be an interval and let $U$ be a unitary operator on $L^2(S)$.

\begin{enumerate}
	\item There exists a Hilbert space \(\mathcal X\subset L^2( S)\), continuously embedded in \(L^2( S)\), such that
	\[
	U\mathcal X=\mathcal X .
	\]
	We also assume that \(\mathcal X\) embeds locally into continuous functions on \(S\), so that
	point evaluations of \(f\in\mathcal X\) and \(Uf\in\mathcal X\) are meaningful.
	Moreover, we assume that \(\mathcal X\) is the form domain of a strictly positive self-adjoint operator \(X\), i.e. $\mathcal X=\mathcal D(X^{1/2})$, and simultaneously its image $\mathcal X=\text{Ran}(X^{-1/2})$, cf.~point \emph{1.}~of Theorem~\ref{thm:mainTheorem}.
	\item The operator \(X\) is represented, in the sense of closed quadratic forms, by\footnote{Since we apply unitary transformations on various domains, we will not distinguish variables using the notation $x$ and $k$ as in the case of the Fourier transform. Henceforth, we will often use the same variable name on both sides.}
	\[
	X=\frac{1}{2}(M_{x^2}+U^*M_{x^2}U),
	\]
	where \(M_{x^2}\) is the self-adjoint multiplication operator
	\(f\mapsto x^2f\). The associated quadratic form is given by
	\[
	 x[f]
	=
	\frac12
	\left(
	\|M_{x^2}^{1/2}f\|_{L^2(S)}^2
	+
	\|M_{x^2}^{1/2}Uf\|_{L^2(S)}^2
	\right),
	\qquad f\in\mathcal X
	\]
	with form domain $
	\mathcal X = D(M_{x^2}^{1/2})\cap U^*D(M_{x^2}^{1/2})$. We assume that this domain is dense in \(L^2(S)\), so that the form is
	densely defined, closed, and non-negative. We define
	\[
	\|f\|_{\mathcal X}^2=2 x[f].
	\]
\end{enumerate}
The assumptions are consistent. If
\(m:S\to\mathbb R\) is measurable, then the multiplication operator
\[
M_m f = mf,
\qquad
\mathcal D(M_m)
=
\{f\in L^2(S): mf\in L^2(S)\},
\]
is self-adjoint on \(L^2(S)\). If \(\widetilde U\) is a unitary operator on \(L^2(S)\), then $\widetilde U^*M_m\widetilde U$ is again self-adjoint, with domain $\mathcal D(\widetilde U^*M_m\widetilde U)
=
\widetilde U^*\mathcal D(M_m)$.
In the case relevant to us, \(m(x)=x^2\), the operator \(M_{x^2}\) is positive. Hence both \(M_{x^2}\) and \(U^*M_{x^2}U\) are positive self-adjoint operators on \(L^2(S)\).

On \(\mathcal X\) we will consider the closed non-negative quadratic forms representing the Hamiltonian's kinetic parts $K_U=U^*M_{x^2}U$, and respectively $K_{U^*}=UM_{x^2}U^*$,
\begin{align*}
	k^U[f]
	=
	\|K_U^{1/2}f\|_{L^2(S)}^2,\qquad
	k^{U^*}[f]
	=
	\|K_{U^*}^{1/2}f\|_{L^2(S)}^2.
\end{align*}

The spectral input is as follows. For every bounded interval
\(I=[a,b]\subset S\), we denote by \( x_I\) the Dirichlet confinement
of the quadratic form \( x\) to \(I\), namely the closed quadratic
form on \(L^2(I)\) obtained by imposing the boundary conditions
\(g(a)=g(b)=0\). The associated positive self-adjoint operator, given by the
first representation theorem, will be denoted by \(X_I\). Similarly,
\( x_I^U\) and \((UXU^*)_I\) denote the corresponding confined form
and operator associated with \(UXU^*\). We denote by \(k_I^U\) and \(k_I^{U^*}\) the Dirichlet
confinements of the kinetic forms \(k^U\) and \(k^{U^*}\), and by
\((K_U)_I\) and \((K_{U^*})_I\) the associated self-adjoint operators.

We restrict attention to settings in
which these confined forms are well defined. We define their spectral floors by
\[
\kappa_I^U
=
\inf_{\substack{g\in\mathcal Q(k_I^U)\\ g\neq 0}}
\frac{k_I^U[g]}{\|g\|_{L^2(I)}^{2}},
\qquad
\kappa_I^{U^*}=\inf_{\substack{g\in\mathcal Q(k_I^{U^*})\\ g\neq 0}}
\frac{k_I^{U^*}[g]}{\|g\|_{L^2(I)}^{2}}.
\]
In other words, $\kappa_I^U=\inf\sigma((K_U)_I)$ and
\(\kappa_I^{U^*}=\inf\sigma((K_{U^*})_I)\). If the confined operators have
compact resolvent, these spectral floors are lowest eigenvalues.

\begin{definition}[Local compatibility with Dirichlet confinement]\label{def:localConfinement}
	Let \( y\) be a closed,
	densely defined, non-negative quadratic form on \(L^2(S)\), with form domain
	\(\mathcal Q( y)\). We say that \( y\) is compatible with
	local Dirichlet confinement if the following property holds:
	
	Whenever \(S\) is decomposed into pairwise disjoint intervals
	\(I_j=(a_j,a_{j+1})\), and \(f\in\mathcal Q( y)\) has a representative
	satisfying \(f(a_j)=0\) at all interior endpoints, then there are respective Dirichlet restrictions \( y_{I_j}\) on \(L^2(I_{j})\) of $ y$ and 
	\(f|_{I_j}\in\mathcal Q( y_{I_j})\) for every \(j\), and
	\[
	 y[f]
	\geq
	\sum_j  y_{I_j}[f|_{I_j}].
	\]
	If equality holds we say that \( y\) is
	additive under Dirichlet confinement.
\end{definition}
The following definition specializes the notation of local compatibility under Dirichlet confinement to the context of uniqueness pairs.

\begin{definition}
	\label{generalCriticalDef}
	Let \((\Lambda,M)\) be a pair of discrete subsets of \(S\) that decompose \(S\)
	and admit strictly increasing enumerations $\Lambda=\{\lambda_j:j\in J\}$, $M=\{\mu_j:j\in J\}$,
	where \(J=\mathbb Z\) in the two-sided case and \(J=\mathbb N\) in the
	one-sided case.
	
	We say that kinetic Dirichlet confinement is admissible for \((\Lambda,M)\)
	if \( k^U\) is compatible with local Dirichlet confinement on the
	decomposition induced by \(\Lambda\), and \( k^{U^*}\) is compatible
	with local Dirichlet confinement on the decomposition induced by \(M\). Thus,
	whenever \(f|_\Lambda=0\) and \(Uf|_M=0\), one has
	\begin{align*}
	 &k^U[f]
	\geq
	\sum_j
	k^U_{[\lambda_j,\lambda_{j+1}]}
	\bigl[f|_{[\lambda_j,\lambda_{j+1}]}\bigr]\quad\text{and}\\
	& k^{U^*}[Uf]
	\geq
	\sum_j k^{U^*}_{[\mu_j,\mu_{j+1}]}
	\bigl[(Uf)|_{[\mu_j,\mu_{j+1}]}\bigr].
	\end{align*}
	We define
	\[
	\rho_j^\Lambda
	=
	\frac{
		\sqrt{\kappa_{[\lambda_j,\lambda_{j+1}]}^U}
	}{
		\max\{|\lambda_j|,|\lambda_{j+1}|\}
	},
	\qquad
	\rho_j^M
	=
	\frac{
		\sqrt{\kappa_{[\mu_j,\mu_{j+1}]}^{U^*}}
	}{
		\max\{|\mu_j|,|\mu_{j+1}|\}
	}.
	\]
	The pair \((\Lambda,M)\) is said to have uniformly supercritical kinetic
	confinement for \(U\) if
	\[
	\inf_j \rho_j^\Lambda>1,
	\qquad
	\inf_j \rho_j^M>1.
	\]
	The corresponding supercritical and subcritical conditions are defined
	analogously, using the appropriate limit superior and limit inferior
	formulations, with the evident one-sided modification when \(J=\mathbb N\).
\end{definition}
\begin{theorem}
	\label{thm:generalizedMainTheorem}
	Let \(U\), \(\mathcal X\), and \(X\) satisfy the assumptions above.
	Assume that kinetic Dirichlet confinement is admissible for \((\Lambda,M)\)
	and that \((\Lambda,M)\) has uniformly supercritical kinetic confinement.
	Then the only function \(f\in\mathcal X\) satisfying
	\(f(\lambda_i)=0\), \(Uf(\mu_i)=0\), \(i\in J\), is \(f=0\).
\end{theorem}

\begin{proof}
	Let \(f\in\mathcal X\) satisfy \(f|_\Lambda=0\) and \(Uf|_M=0\).
	By the uniformly supercritical kinetic confinement assumption, there exists
	\(l>1\) such that
	\[
	\frac{
		\sqrt{\kappa_{[\lambda_j,\lambda_{j+1}]}^{U}}
	}{
		\max\{|\lambda_j|,|\lambda_{j+1}|\}
	}
	\geq l,
	\qquad
	\frac{
		\sqrt{\kappa_{[\mu_j,\mu_{j+1}]}^{U^*}}
	}{
		\max\{|\mu_j|,|\mu_{j+1}|\}
	}
	\geq l
	\]
	for all \(j\). We decompose \(S\) into the intervals determined by
	the zeros of \(f\) and \(Uf\), and follow the lines of the proof of
	Theorem~\ref{thm:mainTheorem}. By admissibility of kinetic Dirichlet
	confinement,
	\begin{align*}
		\norm{f}{\mathcal X}^2
		&=
		\|X^{1/2}f\|_{L^2(S)}^2
		+
		\|(UXU^*)^{1/2}(Uf)\|_{L^2(S)}^2\\
		&\geq
		\frac{1}{2}\sum_{j\in J}\Bigg(
		\|{({K}_U)}_{[\lambda_j,\lambda_{j+1}]}^{1/2}f\|_{L^2([\lambda_j,\lambda_{j+1}])}^2
		+
		\|M_{x^2}^{1/2}f\|_{L^2([\lambda_j,\lambda_{j+1}])}^2\\
		&\hspace{1.5cm}
		+
		\|{(K_{U^*})}_{[\mu_j,\mu_{j+1}]}^{1/2}(Uf)\|_{L^2([\mu_j,\mu_{j+1}])}^2
		+
		\|M_{x^2}^{1/2}Uf\|_{L^2([\mu_j,\mu_{j+1}])}^2
		\Bigg)\\
		&\geq
		\frac{1}{2}\sum_{j\in J}\Bigg(
		\kappa^U_{[\lambda_j,\lambda_{j+1}]}
		\|f\|_{L^2([\lambda_j,\lambda_{j+1}])}^{2}
		+
		\|M_{x^2}^{1/2}f\|_{L^2([\lambda_j,\lambda_{j+1}])}^2\\
		&\hspace{1.5cm}
		+
		\kappa^{U^*}_{[\mu_j,\mu_{j+1}]}
		\|Uf\|_{L^2([\mu_j,\mu_{j+1}])}^{2}
		+
		\|M_{x^2}^{1/2}Uf\|_{L^2([\mu_j,\mu_{j+1}])}^2
		\Bigg)\\
		&\geq
		\frac{1}{2}\sum_{j\in J}\Bigg(
		\left((\rho_j^\Lambda)^2+1\right)
		\|M_{x^2}^{1/2}f\|_{L^2([\lambda_j,\lambda_{j+1}])}^2\\
		&\hspace{3.5cm}
		+
		\left((\rho_j^M)^2+1\right)
		\|M_{x^2}^{1/2}Uf\|_{L^2([\mu_j,\mu_{j+1}])}^2
		\Bigg).
	\end{align*}
	Therefore
	\[
	\begin{aligned}
		\norm{f}{\mathcal X}^2
		&\geq
		\frac{1}{2}(l^2+1)
		\left(
		\|M_{x^2}^{1/2}f\|_{L^2(S)}^2
		+
		\|M_{x^2}^{1/2}Uf\|_{L^2(S)}^2
		\right)\\
		&=
		\frac{1}{2}(l^2+1)\norm{f}{\mathcal X}^2 .
	\end{aligned}
	\]
	Since \(l>1\), we have \(\frac12(l^2+1)>1\), which is impossible unless
	\(\norm{f}{\mathcal X}=0\). Hence \(f=0\).
\end{proof}

Notice that the notation \(L^2(S)\) in Theorem~\ref{thm:generalizedMainTheorem} is understood as the ambient Hilbert space associated with the operator under consideration; for example, for a Sturm-Liouville realization on \(L^2(S,w(x)\d x)\), all confined spectral floors, norms, and multiplication operators are taken with respect to the measure \(w(x)\d x\), see below.

We also record a Hamiltonian version of the criterion, in the spirit of
Lemma~\ref{lem:criticalComparison}. This condition is generally stronger than
the kinetic confinement condition, but it is more natural from the quantum-theoretic point of view. For a bounded interval \(I\subset S\), let \(E_I^X\) and
\(E_I^{UXU^*}\) denote the Dirichlet spectral floors of \(X\) and \(UXU^*\)
on \(I\), respectively:
\[
E_I^X=\inf\sigma(X_I),
\qquad
E_I^{UXU^*}=\inf\sigma((UXU^*)_I).
\]
Equivalently, these are the variational minima of the corresponding confined
quadratic forms. If the confined operators have compact resolvent, then these
spectral floors are their lowest eigenvalues.

We say that \((\Lambda,M)\) has uniformly supercritical Hamiltonian confinement
energy if
\[
\inf_{j\in J}
\frac{\sqrt{E_{[\lambda_j,\lambda_{j+1}]}^X}}
{\max\{|\lambda_j|,|\lambda_{j+1}|\}}
>1\quad\text{and}\quad
\inf_{j\in J}
\frac{\sqrt{E_{[\mu_j,\mu_{j+1}]}^{UXU^*}}}
{\max\{|\mu_j|,|\mu_{j+1}|\}}
>1.
\]

\begin{corollary}[Uniqueness from Hamiltonian confinement energy]
	\label{cor:hamiltonianEnergyUniqueness}
	Let \(U\), \(\mathcal X\), and \(X\) satisfy the assumptions above. Assume
	that kinetic Dirichlet confinement is admissible for \((\Lambda,M)\), and
	that \((\Lambda,M)\) has uniformly supercritical Hamiltonian confinement
	energy. Then the only function \(f\in\mathcal X\) satisfying
	$f(\lambda_i)=0$, $Uf(\mu_i)=0$, $i\in J$ is \(f=0\).
\end{corollary}

\begin{proof}
	As in the proof of Lemma~\ref{lem:criticalComparison}, the bound
	\(x^2\leq \max\{|a|,|b|\}^2\), on \(I=[a,b]\) gives
	\[
	E_I^X\leq \frac12\kappa_I^U+\frac12\max\{|a|,|b|\}^2.
	\]
	Hence \(E_I^X/\max\{|a|,|b|\}^2>1\) implies \(\sqrt{\kappa_I^U}/\max\{|a|,|b|\}>1\), and the same
	argument applied to \(UXU^*\) gives the corresponding estimate for
	\(\kappa_I^{U^*}\). Thus the Hamiltonian confinement condition implies
	uniformly supercritical kinetic confinement, and the result follows from
	Theorem~\ref{thm:generalizedMainTheorem}.
\end{proof}

The abstract results are formulated at the level of local spectral floors. In the Fourier case, the
local spectral floor in Theorem~\ref{thm:generalizedMainTheorem} is the first Dirichlet eigenvalue of the Laplacian and the
resulting estimate is the Wirtinger-Poincar\'{e} inequality. In more
general settings, however, the confined spectral floor is
determined by the corresponding local operator. Thus the Wirtinger-Poincar\'{e} inequality appears as one special instance of the local spectral-floor principle rather than as the
underlying mechanism in all cases. We illustrate this with examples in the following section.

\subsection{Sturm-Liouville-type Operators}
\begin{lemma}[Local confinement for Sturm-Liouville forms] \label{lem:localSLconfinement} Let $S\subseteq\mathbb R$ be an interval and let $w,p:S\to(0,\infty)$ and $q:S\to\mathbb R$ be real coefficient functions. Suppose that $w,p^{-1},q$ are locally $L^1$-integrable inside $S$ and that possible singular endpoints of \(S\) are treated by specific endpoint conditions of the chosen self-adjoint realization. Let $T$ be a positive, self-adjoint operator on $L^2(S,w(x)\d x)$ whose action is locally given by the Sturm-Liouville expression \begin{align}\label{eq:sturmExpression} T f=\frac{1}{w}\left(-(pf')'+qf\right), \qquad f\in\mathcal D(T), \end{align} with endpoint conditions so that $T$ is positive and self-adjoint, and with the Friedrichs condition imposed at singular endpoints. Let \(t\) be the associated closed quadratic form, so that \(\mathcal Q(t)=\mathcal D(T^{1/2})\). Assume that every \(f\in\mathcal Q(t)\) admits a locally continuous representative on \(S\). Assume moreover that the closed form \(t\) is local in the standard Sturm-Liouville form sense: whenever an element of \(\mathcal Q(t)\) vanishes at finitely many interior points of \(S\), its restrictions to the complementary subintervals belong to the corresponding Dirichlet form domains, and the form is additive over such finite decompositions. Let $\{a_j\}\subset S$ be a finite or countable ordered family with no accumulation point in the interior of $S$, and let \(I_j=(a_j,a_{j+1})\). 
	
	Then \(T\) is compatible with Dirichlet confinement on the intervals \(I_j\) in the sense of Definition~\ref{def:localConfinement}. \end{lemma}

Notice that, in this Sturm-Liouville setting, the confined spectral floor is determined by the corresponding local closed form, including its weight, coefficient functions, endpoint conditions, singular endpoints, and boundary behavior.

\begin{proof}[Proof of Lemma~\ref{lem:localSLconfinement}] Let \(f\in\mathcal Q(t)\) vanish at the partition points \(a_j\) lying in the interior of \(S\). Since \(f\) has a locally continuous representative, these pointwise conditions are meaningful. Let \(t_{I_j}\) denote the closed Dirichlet form obtained from \(t\) by restricting the Sturm-Liouville form to \(I_j\), imposing Dirichlet conditions at the finite endpoints of \(I_j\) that lie in the interior of \(S\), and imposing the endpoint condition inherited from \(T\) if \(I_j\) touches an endpoint of \(S\). By the local restriction property of the closed Sturm-Liouville form, \(f|_{I_j}\) belongs to \(\mathcal Q(t_{I_j})\). The Dirichlet conditions at the finite interior endpoints are satisfied because \(f(a_j)=0\), while any endpoint condition at an endpoint of \(S\) is inherited from the original realization. For every finite set \(F\) of indices, locality and additivity over finite decompositions give \[ \sum_{j\in F} t_{I_j}\bigl[f|_{I_j}\bigr] \leq t[f]. \] Since the forms \(t_{I_j}\) are non-negative, the finite partial sums are monotone. Taking the supremum over finite \(F\) yields \[ \sum_j t_{I_j}\bigl[f|_{I_j}\bigr] \leq t[f]. \] This is the compatibility with local Dirichlet confinement required in Definition~\ref{def:localConfinement}. \end{proof}

Our discussion applies immediately to Schr\"{o}dinger operators. The following is a special case of Lemma~\ref{lem:localSLconfinement} with
$w=1$, $p=1$, and $q=V$.
\begin{corollary}[Local confinement for Schr\"{o}dinger forms]
	\label{cor:localSchroedingerConfinement} Under the conditions of Lemma~\ref{lem:localSLconfinement}, let $T$ be a Schr\"{o}dinger operator, whose action is locally
	given by
	\[
	Tf=-f''+Vf,\qquad f\in\mathcal D(T).
	\]
	Then $T$ is compatible with Dirichlet confinement on the intervals $I_j$ in
	the sense of Definition~\ref{def:localConfinement}.
\end{corollary}
\begin{remark}[\citep{RS1}]\label{rem:discreteSpectrumSchroedinger}
	Let $Y=-\frac{\partial^2}{\partial x^2}+V(x)$ be a Schr\"{o}dinger operator on $L^2(\mathbb{R})$, where the potential $V:\mathbb{R}\rightarrow\mathbb{R}\cup\{\infty\}$ is measurable and bounded from below. If $V(x)\rightarrow\infty$ as $\abs{x}\rightarrow\infty$ or $V(x)=\infty$ for $x\notin(a,b)$ (interpreted as a Dirichlet confinement) then $Y$ has purely discrete spectrum.
\end{remark}

It is worth noting that the classical Wirtinger-Poincar\'{e} constant need not give the sharp confinement scale for a weighted Sturm-Liouville operator. In the following example, the preceding confinement criterion is read in the ambient space \(L^2(I,w_\alpha(x)\d x)\).

\begin{example}
Consider on $I=(0,\pi)$ the operator
$$
T_\alpha f=-\frac{1}{w_\alpha(x)}f'' ,
\qquad
w_\alpha(x)=1+\alpha\sin^2(x),
\qquad \alpha>0,$$
with Dirichlet boundary conditions at $0$ and $\pi$. The lowest Dirichlet eigenvalue is characterized by
$$
\lambda_1(\alpha)=\inf_{g\neq 0}
\frac{\int_0^\pi |g'(x)|^2 \d x}
{\int_0^\pi w_\alpha(x)|g(x)|^2\d x}.$$
Testing this quotient with $g(x)=\sin(x)$ gives
$$
\lambda_1(\alpha)
\leq
\frac{\int_0^\pi \cos^2(x)\d x}
{\int_0^\pi (1+\alpha\sin^2(x))\sin^2(x)\d x}
=
\frac{1}{1+\frac34\alpha}
<1.$$

By contrast, the unweighted Wirtinger-Poincar\'{e} constant on $(0,\pi)$ is $1$. Thus a suitable choice of the weight $w$ can lower the actual Dirichlet spectral floor below the value suggested by the unweighted Wirtinger-Poincar\'{e} inequality.
\end{example}
\subsection{The Fractional Fourier Transform}
The fractional Fourier transform $F_\alpha$ generalizes the Fourier transform by providing a representation that corresponds to an intermediate rotation by an angle $\alpha\in[0,2\pi)$ in the time-frequency domain. Numerous works are devoted to the study of the fractional Fourier transform and its applications; see, e.g.,~\citep{YZ,OZK} for a review and references, and~\citep{J} for a discussion of some underlying technical principles.

On the Fourier-symmetric Sobolev space $\mathcal{H}$ the fractional Fourier transform can be defined by the relations
\begin{align*}
	{F_\alpha}[f](k) = e^{-i\alpha H}[f](k),\quad f(x)=e^{i\alpha H}[F_\alpha f](x).
\end{align*}
If $\alpha=0$ no rotation occurs and we obtain the identity map, if $\alpha=\frac{\pi}{2}$ we recover the standard Fourier transform and $\alpha=\pi$ corresponds to spatial reflection: $F_{\pi/2}[f](x)\sim F[f](x)$, $F_\pi[f](x)\sim f(-x)$, both up to a global phase factor. The Hermite functions \(\phi_n\) are eigenvectors of \(F_\alpha\), with
eigenvalues \(e^{-i\alpha(n+1/2)}\). The confinement Hamiltonian associated
with \(F_\alpha\) is
\[
X=X_\alpha=\frac12\left(M_{x^2}+F_\alpha^*M_{x^2}F_\alpha\right),
\]
understood in the sense of closed quadratic forms. Writing $p=-i\frac{\partial}{\partial x}$, $M_{x,\alpha}=F^*_\alpha M_x F_\alpha$ and $p_\alpha=F^*_\alpha p F_\alpha$, we recall the standard commutator relations $[M_x,p]=i$, $[H,M_x]=-ip$ and $[H,p]=iM_x$ and the Heisenberg evolution equations for the conjugated observables
\begin{align*}
	&\frac{\partial}{\partial\alpha}M_{x,\alpha}=i[H,M_{x,\alpha}]=p_\alpha,\\
	&\frac{\partial}{\partial\alpha}p_\alpha=i[H,p_\alpha]=-M_{x,\alpha}.
\end{align*}
Solving this system yields
\begin{align*}
	M_{x,\alpha}&=\cos(\alpha)M_x+\sin(\alpha)p,\\
	p_\alpha&=-\sin(\alpha)M_x+\cos(\alpha)p,
\end{align*}
which justifies that the fractional Fourier transform rotates the phase space and provides a quadratic-form representation for $X$,
\begin{align*}
	2X=(1+\cos^2(\alpha))M_{x^2}+\sin^2(\alpha)p^2+\sin(\alpha)\cos(\alpha)(M_xp+pM_x).
\end{align*}
It is a standard fact in the theory of oscillatory quantum systems that if a quadratic form in $(M_x,p)$ is positive definite, then it can be transformed to a multiple of the Hamiltonian of the standard harmonic oscillator using a symplectic transformation $S:(M_x,p)\mapsto(M_{x'},p')$, $[M_x,p]=[M_{x'},p']$, see~\cite[Chapter 4]{GF} and Williamson's theorem~\citep{JW}. For \(\sin\alpha\neq0\), this implies
\begin{align*}
	2X=\abs{\sin(\alpha)}(M_{x'}^2+p'^2)
\end{align*}
with spectrum
\[
\sigma(X_\alpha)
=
\left\{
|\sin\alpha|\left(n+\frac12\right): n=0,1,2,\ldots
\right\}.
\]
This yields the following uniqueness criterion for the fractional Fourier transform.

\begin{corollary}[Uniqueness pairs of the fractional Fourier transform]
	\label{cor:fracFourierCor}
	Let \(\sin\alpha\neq0\), and let
	\[
	X=\frac12\left(M_{x^2}+F_\alpha^*M_{x^2}F_\alpha\right)
	\]
	in the sense of closed quadratic forms. If one of the following conditions
	holds, then \((\Lambda,M)\) is a uniqueness pair for the fractional Fourier
	transform: the only function \(f\in\mathcal H\) such that $f(\lambda_i)=F_\alpha[f](\mu_i)=0$, $i\in\mathbb Z$
	is \(f=0\).
	
	\begin{enumerate}
		\item The pair \((\Lambda,M)\) has uniformly supercritical
		Hamiltonian confinement energy, that is,
		\[
		\inf_{j\in\mathbb Z}
		\frac{
			\sqrt{E_{[\lambda_j,\lambda_{j+1}]}^X}
		}{
			\max\{|\lambda_j|,|\lambda_{j+1}|\}
		}
		>1
		\]
		and
		\[
		\inf_{j\in\mathbb Z}
		\frac{
			\sqrt{E_{[\mu_j,\mu_{j+1}]}^{F_\alpha X F_\alpha^*}}
		}{
			\max\{|\mu_j|,|\mu_{j+1}|\}
		}
		>1 .
		\]
		
		\item The pair \((\Lambda,M)\) satisfies the following Hamiltonian
		comparison condition:
		\[
		\inf_{j\in\mathbb Z}
		\frac{
			\sqrt{E_{[\lambda_j,\lambda_{j+1}]}^{(0)}}
		}{
			\max\{|\lambda_j|,|\lambda_{j+1}|\}
		}
		>
		\frac{1}{|\sin(\alpha)|}
		\]
		and
		\[
		\inf_{j\in\mathbb Z}
		\frac{
			\sqrt{E_{[\mu_j,\mu_{j+1}]}^{(0)}}
		}{
			\max\{|\mu_j|,|\mu_{j+1}|\}
		}
		>
		\frac{1}{|\sin(\alpha)|},
		\]
		where \(E_{[a,b]}^{(0)}\) denotes the ground-state energy of the ordinary
		confined harmonic oscillator \(H_{[a,b]}\) from~\eqref{confinedSchroed}.
		
		\item The pair \((\Lambda,M)\) satisfies the uniform supercritical density
		conditions
		\[
		\sup_{j\in\mathbb Z}
		\max\{|\lambda_j|,|\lambda_{j+1}|\}
		(\lambda_{j+1}-\lambda_j)
		<
		\pi|\sin\alpha|
		\]
		and
		\[
		\sup_{j\in\mathbb Z}
		\max\{|\mu_j|,|\mu_{j+1}|\}
		(\mu_{j+1}-\mu_j)
		<
		\pi|\sin\alpha|.
		\]
	\end{enumerate}
	On the other hand, if
	\[
	\liminf
	|\lambda_j|(\lambda_{j+1}-\lambda_j)
	>
	\pi|\sin(\alpha)|
	\]
	and
	\[
	\liminf
	|\mu_j|(\mu_{j+1}-\mu_j)
	>
	\pi|\sin(\alpha)|,
	\]
	then \((\Lambda,M)\) is not a uniqueness pair for the fractional Fourier
	transform.
\end{corollary}
\begin{proof} 
	We first prove item~\emph{1.}. We then show that item~\emph{2.}~implies
	item~\emph{3.}. Finally, item~\emph{3.}~is identified with the uniformly
	supercritical kinetic confinement condition for \(F_\alpha\).
	\begin{enumerate}
		\item We first verify that kinetic Dirichlet confinement is admissible. Indeed,
		one may choose the normalized observables
		\begin{align*}
		M_{x'}&=W_\alpha^*M_x W_\alpha=\frac{1}{\sqrt{|\sin(\alpha)|}}M_x,\\
		p'&=W_\alpha^*p W_\alpha=
		\frac{\operatorname{sgn}(\sin(\alpha))}{\sqrt{|\sin(\alpha)|}}
		\bigl(\cos(\alpha)M_x+\sin(\alpha)p\bigr).
		\end{align*}
		The corresponding metaplectic operator may be taken as a dilation followed by
		multiplication by a quadratic phase:
		\[
		(W_\alpha f)(x)
		=
		|\sin(\alpha)|^{1/4}
		\exp\!\left(
		i\frac{\operatorname{sgn}(\sin(\alpha))\cos(\alpha)}{2}x^2
		\right)
		f(\sqrt{|\sin(\alpha)|}\,x).
		\]
		Consequently,
		\begin{align*}
		&f(a)=f(b)=0
		\quad\Longleftrightarrow\\
		&\hspace{1.5cm}(W_\alpha f)\!\left(\frac{a}{\sqrt{|\sin(\alpha)|}}\right)
		=
		(W_\alpha f)\!\left(\frac{b}{\sqrt{|\sin(\alpha)|}}\right)
		=
		0 .
		\end{align*}
		Thus \(W_\alpha\) maps Dirichlet confinement on \([a,b]\) to Dirichlet
		confinement on
		\[
		\left[
		\frac{a}{\sqrt{|\sin(\alpha)|}},
		\frac{b}{\sqrt{|\sin(\alpha)|}}
		\right].
		\]
		
		Now the first assertion follows directly from
		Corollary~\ref{cor:hamiltonianEnergyUniqueness}, applied to \(U=F_\alpha\). Indeed, condition \emph{1.} is the uniformly supercritical Hamiltonian confinement energy condition for \(F_\alpha\).
		\item Item \emph{2.} implies item~\emph{3.}: Assume $\frac{1}{\sin(\alpha)^2}<
		\frac{E_{[a,b]}^{(0)}}{\max\{|a|^2,|b|^2\}}$. As in the proof of Lemma~\ref{lem:criticalComparison} we have that 
		\[
		\frac{1}{\sin(\alpha)^2}<
		\frac{E_{[a,b]}^{(0)}}{\max\{|a|^2,|b|^2\}}
		\leq
		\frac12
		\frac{\nu_{[a,b]}}{\max\{|a|^2,|b|^2\}}
		+\frac12,
		\]
		and therefore
		\[
		\frac{\sqrt{\nu_{[a,b]}}}{\max\{|a|,|b|\}}>\sqrt{\frac{2}{\sin(\alpha)^2}-1}.
		\]
		Since $\nu_{[a,b]}=\pi^2/(b-a)^2$ this implies
		$$\max\{|a|,|b|\}(b-a)<\pi\sqrt{\frac{\sin(\alpha)^2}{2-\sin(\alpha)^2}}\leq\pi |\sin(\alpha)|.$$
		Applying this estimate to
		\([\lambda_j,\lambda_{j+1}]\) and to
		\([\mu_j,\mu_{j+1}]\) yields item \emph{3.}.
		\item It remains to identify the local kinetic spectral floor. The sharp local estimate used for uniqueness comes from the kinetic part operator
		\[
		K_\alpha
		=
		F_\alpha^*M_{x^2}F_\alpha
		=
		\bigl(\cos(\alpha)M_x+\sin(\alpha)p\bigr)^2 .
		\]
		Under the above conjugation, \(K_\alpha\) is mapped to
		\(|\sin(\alpha)|p^2\) on the rescaled interval. Hence the Dirichlet spectral
		floor of \(K_\alpha\) on \([a,b]\) is
		\[
		\kappa_{[a,b]}^\alpha
		=
		|\sin(\alpha)|\,\nu_{[a,b]/\sqrt{|\sin(\alpha)|}}
		=
		\sin^2(\alpha)\frac{\pi^2}{(b-a)^2}.
		\]
		Consequently,
		\[
		\frac{\sqrt{\kappa_{[a,b]}^\alpha}}{\max\{|a|,|b|\}}
		=
		\frac{|\sin(\alpha)|\pi}{\max\{|a|,|b|\}(b-a)}.
		\]
		Condition~\emph{3.} is therefore precisely the uniformly supercritical kinetic
		confinement condition for the fractional Fourier transform. The uniqueness
		claim follows from Theorem~\ref{thm:generalizedMainTheorem}.
	\end{enumerate}
	Finally, the non-uniqueness statement follows from the sharpness result for
	the ordinary Fourier transform in~\citep{KNS}. With the above representation (up to a non-zero constant)
	\[
	F_\alpha f(\xi)
	=
|\sin(\alpha)|^{-1/4}\exp\left(i\frac{\cos(\alpha)\xi^2}{2\sin(\alpha)}\right)\widehat{W_\alpha f}\!\left(\frac{\operatorname{sgn}(\sin(\alpha))\xi}{\sqrt{|\sin(\alpha)|}}\right).
	\]
	Hence
	\[
	f|_\Lambda=0,\qquad F_\alpha f|_M=0
	\]
	is equivalent to
	\[
	(W_\alpha f)|_{\Lambda/\sqrt{|\sin(\alpha)|}}=0,\qquad
	\widehat{W_\alpha f}|_{\operatorname{sgn}(\sin(\alpha)) M/\sqrt{|\sin(\alpha)|}}=0 .
	\]
	If
	\[
	\liminf |\lambda_j|(\lambda_{j+1}-\lambda_j)>\pi|\sin(\alpha)|,
	\qquad
	\liminf |\mu_j|(\mu_{j+1}-\mu_j)>\pi|\sin(\alpha)|,
	\]
	then the rescaled pair $
	\left(\Lambda/\sqrt{|\sin(\alpha)|},
	\operatorname{sgn}(\sin(\alpha)) M/\sqrt{|\sin(\alpha)|}\right)$ is subcritical for the ordinary Fourier transform. By the sharpness result
	of~\cite{KNS}, there exists a nonzero \(h\in\mathcal H\) vanishing on
	\(\Lambda/\sqrt{|\sin(\alpha)|}\) whose Fourier transform vanishes on
	\( \operatorname{sgn}(\sin(\alpha)) M/\sqrt{|\sin(\alpha)|}\). One then chooses $f=W_\alpha^{-1}h\neq0$ and, by construction, $f|_\Lambda=0$ and $F_\alpha f|_M=0$.
\end{proof}
Corollary~\ref{cor:fracFourierCor} shows that the sufficient density threshold
interpolates between the ordinary Fourier threshold at \(\alpha=\pi/2\). At the level of asymptotic density, the limiting condition as
\(\alpha\to0\) or \(\alpha\to\pi\) formally becomes:
$$\limsup\abs{\lambda_j}(\lambda_{j+1}-\lambda_j)=0\ \textnormal{and}\ \limsup\abs{\mu_j}(\mu_{j+1}-\mu_j)=0.$$

\subsection{The Hankel transform}
The Hankel transform of order $\nu>-\frac{1}{2}$, in the unitary convention, is defined as
$$ H_\nu[f](k)=\int_{\mathbb{R}^+}f(r)J_\nu(kr)\sqrt{kr}\d r,\ \textnormal{and}\ f(r)=H_\nu[H_\nu[f]](r),\quad f\in  L^2(\mathbb{R}^+),$$
where $J_\nu$ denotes the Bessel function of first kind of order $\nu$. $H_\nu$ is unitary and involutory, $H_\nu^{-1}=H_\nu$. We consider formally on $L^2(\mathbb{R}^+)$ the \lq\lq{}radial Hamiltonian\rq\rq{}\footnote{This is the one-dimensional radial harmonic oscillator with inverse-square centrifugal term.}
$$X=D_\nu=\frac{1}{2}\left(M_{r^2}+H_\nu M_{r^2}H_\nu\right),$$
defined in the sense of quadratic forms. The respective Hankel-symmetric quadratic form is
\[
 d_\nu[f]
=
\frac12
\left(
\|r f\|_{L^2(\mathbb R^+)}^2
+
\|r H_\nu f\|_{L^2(\mathbb R^+)}^2
\right),
\]
with form domain
\[
\mathcal D_\nu
=
\mathcal D(M_{r^2}^{1/2})
\cap
H_\nu\mathcal D(M_{r^2}^{1/2}).
\]
The operator \(D_\nu\) is the global Hamiltonian associated with the Hankel
transform. Its kinetic part is
\[
L_\nu=H_\nu M_{r^2}H_\nu
=
-\frac{d^2}{dr^2}
+
\frac{\nu^2-\frac14}{r^2}.
\]
Accordingly, the sharp local confinement criterion below is formulated in
terms of the Dirichlet spectral floors of \(L_\nu\), while \(D_\nu\) provides
the global energy space and the Hamiltonian interpretation.
We endow \(\mathcal D_\nu\) with the form norm $\|f\|_{\mathcal D_\nu}^2
=
2 d_\nu[f]$.
\begin{proposition}\label{radialHamiltonian}
	The following assertions hold for the Hankel-symmetric space $\mathcal{D}_\nu$ with $\nu>0$:
	\begin{enumerate}
		\item $\cD_\nu\subset L^2(\mathbb{R}^+)$ and $\cD_\nu$ is the form domain of the self-adjoint and positive operator 
		$$D_\nu=\frac{1}{2}(L_\nu+M_{r^2})\quad\textnormal{with}\quad L_\nu=-\left(\frac{\partial^2}{\partial r^2}-\frac{\nu^2-\frac{1}{4}}{r^2}\right),$$
		where \(L_\nu\) is understood as the Friedrichs realization on
		\(L^2(\mathbb R^+)\).
		\item $D_\nu$ has eigenvectors that are Laguerre functions of the form $$\varphi_{n,\nu}=r^{\nu+\frac{1}{2}}e^{-\frac{r^2}{2}}L_{n,\nu}(r^2),\ n\geq0,$$
		where $L_{n,\nu}$ denotes the associated Laguerre polynomial.
		The corresponding eigenvalues are $E_{n,\nu}=(2n+\nu+1)$.
		\item The Laguerre functions form an orthogonal basis of \(L^2(\mathbb R^+)\)
		and an orthogonal basis of \(\mathcal D_\nu\) with respect to the form inner
		product.
		\item $D_\nu$ is self-dual under the Hankel transform: ${H_\nu}D_\nu{H_\nu}^{-1}=D_\nu$.
		
		\item For all \(f\in\mathcal D_\nu\), \(f\neq0\), it holds that
		\[
		\|f\|_{\mathcal D_\nu}
		\geq
		\sqrt{2(\nu+1)}\,\|f\|_{L^2(\mathbb R^+)}.
		\]
		Equality is achieved if and only if \(f\) is a multiple of
		\(\varphi_{0,\nu}\).
	\end{enumerate}
\end{proposition}
These facts are known but as notational conventions vary, we provide short proofs of \emph{(1)}, \emph{(2)}, \emph{(3)} for the reader's convenience. Property \emph{(4)} follows from the definition of \(D_\nu\), and
\emph{(5)} follows from the spectral lower bound \(E_{0,\nu}=\nu+1\). Refer to~\citep{GW} for details.
\begin{proof} Recall that the Bessel functions of the first kind solve the Bessel differential equation $$\frac{\partial^2J_\nu}{\partial r^2}+\frac{1}{r}\frac{\partial J_\nu}{\partial r}+\left(1-\frac{\nu^2}{r^2}\right)J_\nu=0,\quad x>0.$$
	\begin{enumerate}
		\item Substituting $x=rk$ and applying the product rule, Bessel's equation transforms to an eigenvalue equation
		$$\left[-\frac{\partial^2}{\partial r^2}+\frac{\nu^2-\frac{1}{4}}{r^2}\right]\left(\sqrt{kr}J_\nu(kr)\right)=k^2\left(\sqrt{kr}J_\nu(kr)\right).$$
		%
		%transforms the Bessel equation to
		%$$r^2J_\nu(kr)=-\left(\frac{\partial^2}{\partial k^2}+\frac{1}{k}\frac{\partial}{\partial k}-\frac{\nu^2}{k^2}\right)J_\nu(kr),$$
		%which yields via a short computation to
		Applying the Hankel transform to this identity yields
		$$H_\nu[r^2f](k)=-\left(\frac{\partial^2}{\partial k^2}-\frac{\nu^2-\frac{1}{4}}{k^2}\right)H_\nu[f](k).$$
		Thus $H_\nu M_{r^2}H_\nu=L_\nu$ as claimed.
		\item With ladder operators chosen as
		\begin{align*}
			a_\nu=\frac{\partial}{\partial r}+r-\frac{\nu+1/2}{r},\quad  a_\nu^\dagger=-\frac{\partial}{\partial r}+r-\frac{\nu+1/2}{r},
		\end{align*}
		direct computation shows $D_\nu=\frac{1}{2}\left(a_\nu^\dagger a_\nu+2\nu+2\right)$. Making use of the Laguerre identities
		\begin{align*}
			\frac{\partial}{\partial x}L_{n,\nu}(x)&=-L_{n-1,\nu+1}(x),\\
			L_{n,\nu}(x)&=L_{n,\nu+1}(x)-L_{n-1,\nu+1}(x),\\
			nL_{n,\nu}(x)&=(n+\nu)L_{n-1,\nu}(x)-xL_{n-1,\nu+1},
		\end{align*}
		one verifies the action of ladder operators
		\begin{align*}
			a_\nu\varphi_{n,\nu}=-2\varphi_{n-1,\nu+1},\quad
			a^\dagger_\nu\varphi_{n-1,\nu+1}=-2n\varphi_{n,\nu}.
		\end{align*}
		This yields the claimed eigenvector-eigenvalue relationship for $n,\nu$.
		\item This is a consequence of the orthogonality of the associated Laguerre polynomials $L_{n,\nu}(x)$ with respect to the weighting function $x^\nu e^{-x}$.
	\end{enumerate}
\end{proof}

We next estimate the confined spectral floors of the Bessel kinetic operator
\(L_\nu\) restrict to \(\nu>0\), so that
the Bessel operator is taken with its Friedrichs realization at the singular
endpoint. In the one-sided Hankel setting we set
\(\lambda_0=\mu_0=0\) and interpret the confinement on the first intervals
\([0,\lambda_1]\) and \([0,\mu_1]\) with the Friedrichs boundary condition at
the singular endpoint \(0\). The zeros at positive points impose ordinary
Dirichlet walls.

This leads to a sufficient condition for uniqueness pairs for the Hankel transform
by Theorem~\ref{thm:generalizedMainTheorem} and Corollary~\ref{cor:hamiltonianEnergyUniqueness}. In the sufficient condition below, the
estimates are applied to intervals bounded away from \(0\). The interval adjacent to $0$ is treated separately through the Friedrichs confined spectral floor for $\nu>0$. As discussed for the Fourier transform, for compact subintervals of \((0,\infty)\), elements of \(\mathcal D_\nu\) have locally continuous representatives and point evaluations at positive points are well
defined.
\begin{corollary}[Uniqueness pairs of the Hankel transform]
	\label{cor:hankelCor}
	Let \(\nu>0\), and let
	\[
	\Lambda=\{0=\lambda_0<\lambda_1<\lambda_2<\cdots\},
	\qquad
	M=\{0=\mu_0<\mu_1<\mu_2<\cdots\}
	\]
	be two increasing sequences that decompose \([0,\infty)\). If one of the following conditions
	holds, then \((\Lambda,M)\) is a uniqueness pair for the Hankel
	transform: the only
	function \(f\in\mathcal D_\nu\) satisfying $f(\lambda_j)=0$, $H_\nu f(\mu_j)=0$ for all $j\geq1$ is \(f=0\).
	\begin{enumerate}
		\item The pair \((\Lambda,M)\) has uniformly supercritical
		kinetic Dirichlet confinement, that is, for the Dirichlet spectral floor \(\kappa_I^{L_\nu}\) of the
		Friedrichs realization of \(L_\nu\) on \(I\), it holds that
		\[
		\inf_{j\geq0}
		\frac{\sqrt{\kappa_{[\lambda_j,\lambda_{j+1}]}^{L_\nu}}}
		{\lambda_{j+1}}
		>1,
		\qquad
		\inf_{j\geq0}
		\frac{\sqrt{\kappa_{[\mu_j,\mu_{j+1}]}^{L_\nu}}}
		{\mu_{j+1}}
		>1.
		\]
			In particular, this conclusion holds if the first intervals satisfy
		\[
		\frac{\sqrt{\kappa_{[0,\lambda_1]}^{L_\nu}}}{\lambda_1}>1,
		\qquad
		\frac{\sqrt{\kappa_{[0,\mu_1]}^{L_\nu}}}{\mu_1}>1,
		\]
		and if, for the remaining intervals,
		\begin{align*}
		\inf_{j\geq1}
		\frac{1}{\lambda_{j+1}^{2}}
		\left[
		\frac{\pi^{2}}{(\lambda_{j+1}-\lambda_j)^{2}}
		+
		\inf_{r\in[\lambda_j,\lambda_{j+1}]}
		\frac{\nu^{2}-\tfrac14}{r^{2}}
		\right]
		>1,\\
		\inf_{j\geq1}
		\frac{1}{\mu_{j+1}^{2}}
		\left[
		\frac{\pi^{2}}{(\mu_{j+1}-\mu_j)^{2}}
		+
		\inf_{r\in[\mu_j,\mu_{j+1}]}
		\frac{\nu^{2}-\tfrac14}{r^{2}}
		\right]
		>1 .
		\end{align*}
		\item The pair \((\Lambda,M)\) has uniformly supercritical
		Hamiltonian confinement energy, that is, for the Dirichlet spectral floor \(E_I^{D_\nu}\) of the
		Friedrichs realization of \(D_\nu\) on \(I\), it holds that
		\[
		\inf_{j\geq0}
		\frac{\sqrt{E_{[\lambda_j,\lambda_{j+1}]}^{D_\nu}}}
		{\lambda_{j+1}}
		>1,
		\qquad
		\inf_{j\geq0}
		\frac{\sqrt{E_{[\mu_j,\mu_{j+1}]}^{D_\nu}}}
		{\mu_{j+1}}
		>1.
		\]
			In particular, this conclusion holds if the first intervals satisfy
		\[
		\frac{\sqrt{E_{[0,\lambda_1]}^{D_\nu}}}{\lambda_1}>1,
		\qquad
		\frac{\sqrt{E_{[0,\mu_1]}^{D_\nu}}}{\mu_1}>1,
		\]
			and if, for the remaining intervals,
		\begin{align*}
		\inf_{j\geq1}
		\frac{1}{\lambda_{j+1}^{2}}
		\left[
		\frac{\pi^{2}}{2(\lambda_{j+1}-\lambda_j)^{2}}
		+
		\frac12
		\inf_{r\in[\lambda_j,\lambda_{j+1}]}
		\left(
		\frac{\nu^{2}-\tfrac14}{r^{2}}+r^{2}
		\right)
		\right]
		>1,\\
		\inf_{j\geq1}
		\frac{1}{\mu_{j+1}^{2}}
		\left[
		\frac{\pi^{2}}{2(\mu_{j+1}-\mu_j)^{2}}
		+
		\frac12
		\inf_{r\in[\mu_j,\mu_{j+1}]}
		\left(
		\frac{\nu^{2}-\tfrac14}{r^{2}}+r^{2}
		\right)
		\right]
		>1 .
		\end{align*}
	\end{enumerate}
\end{corollary}

\begin{proof} The admissibility of the corresponding kinetic confinement follows from
	Lemma~\ref{lem:localSLconfinement}, applied to the Bessel
	Sturm-Liouville expression with the Friedrichs condition at \(0\).
	\begin{enumerate}
		\item 
	The kinetic confinement operator associated with \(H_\nu\) is
	\(L_\nu\), both in the original and transformed representations. The stated
	spectral-floor condition is therefore the uniformly
	supercritical kinetic confinement condition of
	Theorem~\ref{thm:generalizedMainTheorem}. The final sufficient condition
	follows from the lower bound
	\[
	\kappa_{[a,b]}^{L_\nu}
	\geq
	\frac{\pi^2}{(b-a)^2}
	+
	\inf_{r\in[a,b]}
	\frac{\nu^2-\frac14}{r^2}
	\]
	on intervals compactly contained in \((0,\infty)\).
	\item The first assertion follows from
	Corollary~\ref{cor:hamiltonianEnergyUniqueness}, applied with
	\(U=H_\nu\) and \(X=D_\nu\). The displayed sufficient condition follows from
	the variational lower bound
	\[
	E_{[a,b]}^{D_\nu}
	\geq
	\frac{\pi^2}{2(b-a)^2}
	+
	\frac12
	\inf_{r\in[a,b]}
	\left(
	\frac{\nu^2-\tfrac14}{r^2}+r^2
	\right)
	\]
	on intervals compactly contained in \((0,\infty)\), together with the
	corresponding Friedrichs spectral floors on the first intervals adjacent to
	\(0\).
		\end{enumerate}
\end{proof}
For intervals compactly contained in \((0,\infty)\), these spectral floors are
bounded from below by the ordinary Dirichlet Laplacian floor together with the
inverse-square Bessel potential. Thus the resulting estimate is a
Sturm-Liouville refinement of the Wirtinger-Poincar\'e bound.

The intervals adjacent to the singular endpoint \(0\) are different. There the
Friedrichs realization of \(L_\nu\) encodes the admissible boundary behavior at
the inverse-square singularity, and the corresponding spectral floors reflect
the Hardy-type control associated with the Bessel operator at \(0\).

\section{Conclusion}

This work develops a spectral framework for uniqueness pairs associated with unitary transformations. Zeros of a function and of its transform are interpreted as Dirichlet-type confinement in two conjugate representations, converting two-sided sampling problems into spectral lower-bound problems for confined Hamiltonians.

In the Fourier case, this spectral viewpoint recovers the known supercritical uniqueness criterion through a variational argument: the relevant local spectral floor is precisely the Dirichlet Laplacian constant from the Wirtinger-Poincar\'{e} inequality. For more general transforms, the same mechanism is governed by the full confined operator associated with the transform, including any potential terms, weights, or endpoint effects.

This distinction becomes visible in the extensions considered here. For the
fractional Fourier transform, the criterion reflects the rotation of phase
space and leads to a rescaling of the critical threshold by the corresponding
rotation parameter. For the Hankel transform, the relevant operator is a
radial harmonic oscillator with an inverse-square Bessel potential and a
Friedrichs boundary condition at the singular endpoint. The spectral floors in
this case contain not only a Poincar\'e-type contribution from ordinary
Dirichlet confinement, but also the Hardy-type endpoint behavior and the
Sturm-Liouville structure of the radial operator.

The article connects sampling theory, harmonic analysis, and mathematical physics by interpreting two-sided uniqueness as a spectral confinement problem. Sampling criteria become lower bounds on confined spectral floors, while the physical interpretation is a confinement form of the uncertainty principle: a nonzero state cannot be excessively confined in two conjugate representations.

\section*{Acknowledgements}

This work was carried out during a research stay hosted by Andrii Bondarenko
and Kristian Seip. The author is grateful to them for valuable discussions,
and also thanks the Department of Mathematical Sciences at NTNU for its warm
hospitality and supportive environment. This work was supported by the Swiss
National Science Foundation, SNSF grant no.~CRSK-2\_229036.

\bibliographystyle{plainnat}
\bibliography{sn-bibliography}

\end{document}